\documentclass{article}

\usepackage{amsmath}
\usepackage{graphicx}
\usepackage{enumitem}
\usepackage{caption,subcaption}
\usepackage{float}
\usepackage{calc}
\usepackage[calc,showdow,en-US]{datetime2}
\usepackage{relsize}
\usepackage{layout}
\usepackage{geometry}

\DTMsetstyle{en-US}

\usepackage{relsize}

\newtheorem{theorem}{Theorem}

\title{A Mixed Graph Achieving a Moore-like Bound}

\author{Geoffrey Exoo \\
        Department of Mathematics and Computer Science \\
        Indiana State University \\
        Terre Haute, IN 47802 \\
        {\tt ge@cs.indstate.edu}}

\begin{document}

\large

\maketitle

\begin{center}
\end{center}

\thispagestyle{empty}
\pagestyle{empty}
\begin{abstract}
Mixed graphs have both directed and undirected edges.  A mixed cage
is a regular mixed graph of given girth with minimum possible order.
In this paper we construct a mixed cage of order $30$ that achieves the mixed
graph analogue of the Moore bound for degree 3, out-degree 1, and
girth 6.
\end{abstract}

\section{Notation and Terminology}

A mixed graph is a graph with both directed and undirected edges.
We refer to directed
edges as {\em arcs} and undirected edges as {\em edges}.
The {\em degree} of a vertex $v$ in a
mixed graph $G$ is the number of edges incident with $v$, whereas
the {\em in-degree} and {\em out-degree} are the numbers of arcs
incident to and from $v$.  $G$ is {\em regular} if all three degrees
are constant as $v$ ranges over $V(G)$.

A {\em cycle} in a mixed graph is a sequence of vertices $v_0, v_1, \cdots, v_k$ such that
there are no repeated vertices except that $v_0 = v_k$,
each pair of consecutive vertices $(v_i,v_{i+1})$ is
either an edge or an arc, and there are no repeated edges or arcs.
The {\em girth} of a mixed graph is the length of a shortest cycle.  Note that this definition
considers the possibility of $1$-cycles (loops) and $2$-cycles.

Mixed graphs have been studied in the context of the degree/diameter
problem.  We follow the notation used in that literature and denote
the degree and outdegree of a regular mixed graph by $r$ and $z$, respectively.

An $(r,z,g)$-graph
is a regular mixed graph with degree $r$, out-degree $z$, and girth $g$.
An $(r,z,g)$-cage is an $(r,z,g)$-graph of minimum possible order.
We denote this minimum order by $f(r,z,g)$.

\medskip

Recall \cite{ddmsurvey} that the Moore bound for an $r$-regular graph of diameter $d$ is given by:

\begin{equation}
n(r,d) \ge \frac{r(r-1)^{d}-2}{r-2}
\end{equation}

In \cite{ahm} 
Araujo-Pardo, Hernández-Cruz, and Montellano-Ballesteros consider
the problem of finding mixed cages.
They focus on the case $z=1$ and determine a lower bound for $f(r,1,g)$
based on the Moore bound.
Their idea is to
attach undirected Moore trees to each vertex of a directed path of length $g{-}1$,
choosing trees whose depth is as large as possible while still guaranteeing that
all tree vertices are distinct.

So let $v_0, v_1, \cdots , v_{g-1}$ be the vertices of a directed path of length $g{-}1$.
Using edges, attach a Moore tree of depth $i$ to both $v_i$ and $v_{g-1-i}$
for $0 \leq i \leq \lfloor g/2 \rfloor$.  Note: if $i = g{-}i{-}1$ we attach only
one tree.
The base path contains arcs and the Moore trees contain edges.
This gives the following bound.

\begin{theorem}[The AHM Bound \cite{ahm}]

\begin{equation*}
f(r,1,g) \geq \mathlarger{\mathlarger{\sum}}_{i=0}^{g-1}n(r,\, {\tt min}(i,g{-}i{-}1))
\end{equation*}

\end{theorem}

\medskip

\begin{theorem}
\begin{equation*}
f(3,1,6)  = 30
\end{equation*}
\end{theorem}
{\bf Proof.}
The AHM bound for $r=3$, $z=1$, and $g=6$ is $30$.  The graph $G$ shown
in the figure has order $30$ and has the require parameters.

\medskip
The graph can be describe algebraically as follows.
Let
\begin{equation*}
V(G) = \{ v(i,j) \,|\, 0 {\leq} i {<} 3,\, 0 {\leq} j {<} 10 \}
\end{equation*}
be the vertex set.
For $0 \leq i < 3$, let
\begin{equation*}
V_i = \{ v(i,j) \,|\, 0 {\leq} j {<} 10\}. 
\end{equation*}

Each $V_i$ induces a directed $10$-cycle.  These three cycles
are shown in figure.  Let $V_0$ be the set of vertices
in the lower partion of the figure and let $V_1$ and $V_2$ be the sets of
vertices on the outer and inner (resp.) cycles in the upper part of the figure.
Imagine the vertices labeled such that
vertex $v_{i,0}$ is the rightmost vertex on each of the three directed $10$-cycles
and that the other vertices
are labeled in counter-clockwise order.
The arcs and edges are given as follows,
where second indices are computed module $10$:

\begin{enumerate}[label=(\alph*) ]
\item arc($v(0,j) ,\, v(0,j{+}1)$)
\item edge($v(0,j) ,\, v(1,j)$)
\item edge($v(0,j) ,\, v(2,j{+}5)$)
\item edge($v(1,j) ,\, v(2,j{+}2)$)
\item edge($v(1,j) ,\, v(2,j{-}2)$)
\end{enumerate}

In the figure, edges of type (b) and (c) are not drawn but indicated by color:
vertices in $V_0$ are adjacent to those vertices in $V_1$ and $V_2$ that have
the matching color.

\noindent
{\bf Note:}
Recently another graph achieving the AHM bound was found by the author.
This graph is for the case $(6,1,6)$.
The graph has order $90$, which is the AHM bound, and can be constructed as a lift of the
complement of the line graph of $K_6$.  The adjacency matrix for the graph can be found
at the following location.

\begin{center}
\verb+https://cs.indstate.edu/ge/MixedCages/g90.txt+
\end{center}

\begin{figure}
\centering
\includegraphics[scale=0.5]{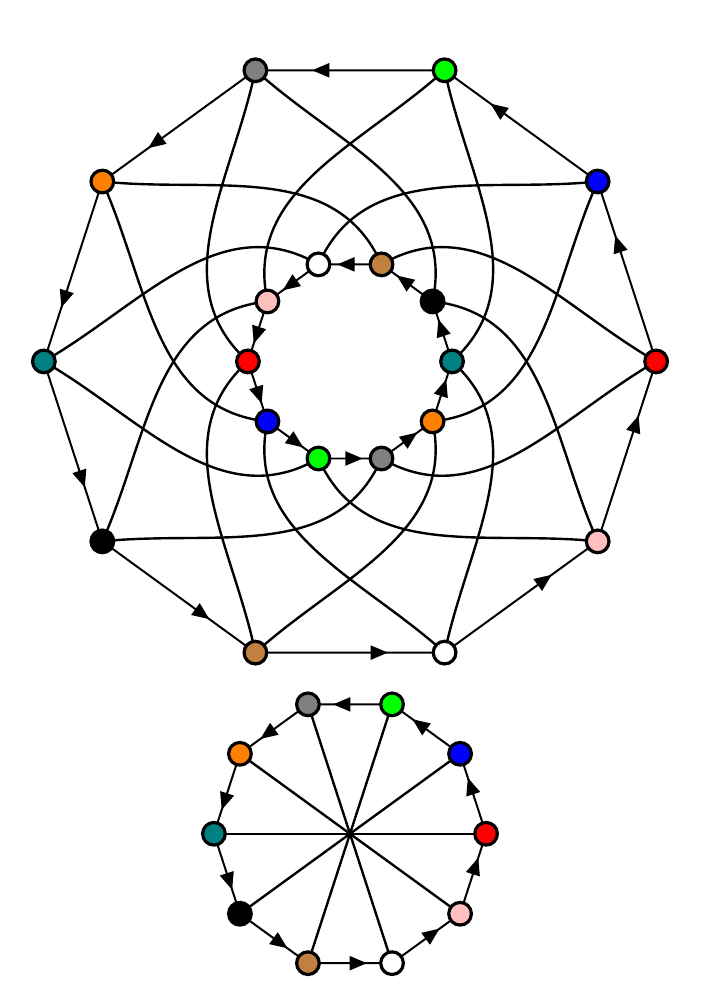}
\caption{The unique smallest $(3,1,6)$-graph of order $30$.  Vertices in the lower figure are
adjacent to vertices in the upper figure that have the matching color.  The automorphism
group has two generators: a rotation of $\pi/5$ (in both figures) and an involution that
transposes the inner and outer cycles in the upper figure.
These generators commute, so the group is $Z_2 \times Z_{10}$.
}
\label{g316fig}
\end{figure}

\end{document}